\documentclass[11pt,reqno]{amsart}
\usepackage{amssymb,latexsym}
\usepackage{graphicx}
\usepackage{fancyhdr,amssymb}
\usepackage{url}		
\usepackage{graphics}
\usepackage[hang]{subfigure}
\usepackage{pgf,tikz}

\def \no {\noindent}
\def \beq {\begin{equation}}
\def \eeq {\end{equation}}
\theoremstyle{plain}
\numberwithin{equation}{section}
\newtheorem{thm}{Theorem}[section]
\newtheorem{theorem}[thm]{Theorem}
\newtheorem{lemma}[thm]{Lemma}

\newtheorem{definition}[thm]{Definition}

\begin{document}
\fancyhead{}
\renewcommand{\headrulewidth}{0pt}
\fancyfoot{}
\fancyfoot[LE,RO]{\medskip \thepage}
\fancyfoot[LO]{\medskip MONTH YEAR}
\fancyfoot[RE]{\medskip VOLUME , NUMBER }

\setcounter{page}{1}

\title[Closed-form representation of Fibonacci numbers]{A note on closed-form representation of Fibonacci numbers using Fibonacci trees}
\author{Indhumathi Raman}
\address{School of Information Technology and Engineering\\ VIT University\\ Vellore, India}
\email{indhumathi.r@vit.ac.in}

\begin{abstract}
In this paper, we give a new representation of the Fibonacci numbers. This is achieved using Fibonacci trees. With the help of this representation, the $n^{th}$ Fibonacci number can be calculated without having any knowledge about the previous Fibonacci numbers.
\end{abstract}

\maketitle

\date{\empty}

\section{Introduction}

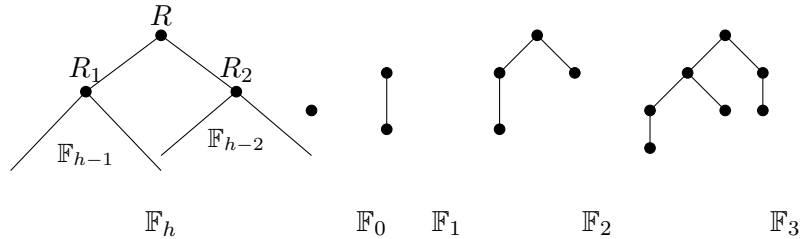
\begin{figure}[ht]
\begin{center}
\begin{tikzpicture}
\draw   (0,3.25) -- (-1,4) -- (-2,3.25)
    (-1,2.2) -- (-2,3.25) -- (-3,2.2)
    (1,2.4) -- (0,3.25) -- (-1,2.4);

\draw   (-2,2.4) node {\small $\mathbb{F}_{h-1}$}
    (0,2.6) node {\small $\mathbb{F}_{h-2}$}
    (-1,1.5) node {$\mathbb{F}_h$};

\filldraw[black,draw=black]
    (-2,3.25) node[above] {$R_1$}  circle (2pt)
    (0,3.25) node[above] {$R_2$}  circle (2pt)
    (-1,4) node[above] {$R$} circle (2pt);

\hspace*{0.8cm}

\filldraw[black, draw=black]
    (1,3) circle (2pt);

\filldraw[black, draw=black]
    (2,3.5) circle (2pt)
    (2,2.75) circle (2pt);
\draw   (2,2.75) -- (2,3.5);

\filldraw[black, draw=black]
    (4,4) circle (2pt)
    (3.5,3.5) circle (2pt)
    (4.5,3.5) circle (2pt)
    (3.5,2.75) circle (2pt);
\draw   (3.5,2.75) -- (3.5,3.5) -- (4,4) -- (4.5,3.5);

\filldraw[black, draw=black]
    (6.5,4) circle (2pt)
    (6,3.5) circle (2pt)
    (7,3.5) circle (2pt)
    (5.5,3) circle (2pt)
    (6.5,3) circle (2pt)
    (7,3) circle (2pt)
    (5.5,2.5) circle (2pt);

\draw   (5.5,2.5) -- (5.5,3) -- (6,3.5) -- (6.5,3)
        (6,3.5) -- (6.5,4) -- (7,3.5) -- (7,3);

%
%

\draw   (1,1.5) node {$\mathbb{F}_0$}
        (2,1.5) node {$\mathbb{F}_1$}
        (4,1.5) node {$\mathbb{F}_2$}
        (6.5,1.5) node {$\mathbb{F}_3$};
\end{tikzpicture}
\caption{Recursive construction and examples of Fibonacci Trees} \label{FTexmpls}
\end{center}
\end{figure}

A Fibonacci tree is a rooted binary tree in which for
every non-leaf vertex $v$, the heights of the subtrees, rooted at
the left and right child of $v$, differ by exactly one.
A formal recursive definition of the Fibonacci tree (denoted by $\mathbb{F}_h$ if its height is $h$) is given below.
\begin{definition} \label{Fibrecurs}
$\mathbb{F}_0 := K_1, ~\mathbb{F}_1 := K_2$. For $h \geq 2$, $\mathbb{F}_h$ is obtained by taking a copy of $\mathbb{F}_{h-1}$, a copy of $\mathbb{F}_{h-2}$, a new vertex $R$ and joining $R$ to the roots of $\mathbb{F}_{h-1}$ and $\mathbb{F}_{h-2}$.
\end{definition}
\no Figure \ref{FTexmpls} shows this construction and a few small Fibonacci trees.

The above recursive definition implies that the number of vertices in $\mathbb{F}_h$ is
$|V(\mathbb{F}_h)| = |V(\mathbb{F}_{h-1})| + |V(\mathbb{F}_{h-2})| + 1$. On solving this recurrence relation, we get
$|V(\mathbb{F}_h)| = f(h+2) - 1$, where $f(i)$ is the $i^{th}$ number in the Fibonacci sequence,
$f(0) = 1, f(1) = 1, ~f(n) = f(n-1) + f(n-2)$; this justifies the terminology Fibonacci tree.
The Fibonacci tree is the one with minimum number of vertices among the class of AVL trees (\cite{avl}). Several properties of Fibonacci trees have been investigated: for eg., Fibonacci numbers of Fibonacci trees has been studied in \cite{FibFib}, Optimality of Fibonacci numbers is discussed in \cite{Optim}, asymptotic properties of Balaban's index for Fibonacci trees has been explored in \cite{Balaban}, Zeckendorf representation of integers is given in \cite{zeck}. In this short paper, we represent the number of vertices of $\mathbb{F}_h$ in {\it closed-form}\footnote{A closed-form is one which gives the value of a sequence at index $n$ in terms of only one parameter, $n$ itself.} by observing the number of vertices at each level of $\mathbb{F}_h$. Such a calculation helps us to give a closed-form representation of $n^{th}$ Fibonacci number for every $n \geq 2$.

\section{Closed-form representation of Fibonacci numbers}
There are several closed-form representations of the Fibonacci numbers. We state a few below.

\begin{itemize}

\item $f(n) = \frac{(1+\sqrt{5})^n - (1-\sqrt{5})^n}{2^n \sqrt{5}}$

It was derived by Binet \cite{Binet} in 1843, although the result was known to Euler, Daniel Bernoulli, and de Moivre more than a century earlier.

\item $B(x) = \displaystyle \sum_{k=0}^{\infty} b_kx^k$

In the above generating function for the Fibonacci numbers' the value of $b_k$ gives the $k^{th}$ Fibonacci number. However expanding the generating function involves tedious calculations.

\item $f_n = \text{round}\left(\frac{5+\sqrt{5}}{10}{\left(\frac{1+\sqrt{5}}{2}\right)}^n\right)$

It was also derived by Binet \cite{Binet} where the function $round()$ rounds the simplified expression up or down to an integer.
\end{itemize}

In this section, we give a simpler closed-form combinatorial representation for Fibonacci numbers. To do so, we first give a closed-form representation for the number of vertices $|V(\mathbb{F}_h)|$ of $\mathbb{F}_h$ (the Fibonacci tree of height $h$). 
One of the most powerful methods of solving equations is called Guess-and-Check - We can use any
method to figure out a solution to an equation and then determine whether it is right by substituting it
in the equation and checking the equation is solved by the guessed solution. This method is particularly
powerful when used with differential equations, recurrence equations and simultaneous equations -
Such a method is employed in the following lemma to calculate the number of vertices in a particular level of $\mathbb{F}_h$ and thereafter we sum the number of vertices over the levels to get $|V(\mathbb{F}_h)|$.

\begin{lemma} \label{formulalemma}
Let $\mathbb{F}_h$ be a Fibonacci tree of height $h$ and let $k$ be an integer such that $0 \leq k \leq h$. The number of vertices $N(h,k)$ at level $k$ of $\mathbb{F}_h$ is given by $$ N(h,k) = \displaystyle \sum_{i=0}^{h-k} \binom{k}{h-k-i}$$
\end{lemma}

\no {\bf Proof}: We prove the lemma by induction on $k$. For $k=0$ we have $N(h,0) = \displaystyle \sum_{i=0}^{h} \binom{0}{h-i}$. Using the convention $\binom{n}{r} = 0 ~\text{if}~ n < r$, we have $N(h,0)= \binom{0}{0} = 1$. This is true since the root of $\mathbb{F}_h$ is the only vertex at level $0$. Further proceeding, from the recursive definition of $\mathbb{F}_h$, we have
\begin{eqnarray*}
N(h,k) &= & N(h-1,k-1) + N(h-2,k-1) \\
&= & \displaystyle \sum_{i=0}^{h-k} \binom{k-1}{h-k-i} + \displaystyle \sum_{j=0}^{h-k-1} \binom{k-1}{h-k-j-1} \\
&= & \displaystyle \sum_{i=0}^{h-k} \binom{k-1}{h-k-i} + \displaystyle \sum_{j=0}^{h-k} \binom{k-1}{h-k-j-1} - \binom{k-1}{-1} \\
&= & \displaystyle \sum_{i=0}^{h-k} \left( \binom{k-1}{h-k-i} +  \binom{k-1}{h-k-i-1}\right) ~\text{since} \binom{n}{r} = 0 ~\text{if}~ r < 0 \\
&= & \displaystyle \sum_{i=0}^{h-k} \binom{k}{h-k-i} \\
\end{eqnarray*}
In Step $3$ of the above equation, we add and subtract $\binom{k-1}{h-k-j-1}$ for $j = h-k$. This proves the lemma.~~ \hfill{$\Box$}

\vspace*{0.5cm}
\no The number of vertices in any tree is the sum of the vertices at its levels. In particular, $|V(\mathbb{F}_h)| = \displaystyle \sum_{k=0}^{h} N(h,k)$. Hence we have the following lemma.

 \begin{lemma} \label{formulavertices}
Let $\mathbb{F}_h$ be the Fibonacci tree of height $h$, then the number of vertices $|V(\mathbb{F}_h)|$ of $\mathbb{F}_h$ is $\displaystyle \sum_{k=0}^{h} \sum_{i=0}^{h-k}  \binom{k}{h-k-i}$.   \hfill{$\Box$}
\end{lemma}

The above theorem helps us to derive a closed-form representation for the Fibonacci numbers. This representation is in contrast to the recurrence relation form, which has certain previous values of the sequence as parameters. We know that $|V(\mathbb{F}_h)| = f(h+2)-1$. Equivalently $f(n) = 1+ |V(\mathbb{F}_{n-2})|$.

\begin{theorem} \label{formulaFib}
Let $f(n)$ be the $n^{th}$ number in the Fibonacci sequence starting with $f(0)=1$ and $f(1)=1$. Then for $n \geq 2$, $$f(n) = 1+ \displaystyle \sum_{k=0}^{n-2} \sum_{i=0}^{n-k-2} \binom{k}{n-k-i-2}$$
\end{theorem}

\no {\bf Proof}: An immediate consequence of Lemma \ref{formulavertices}. \hfill{$\Box$}

\noindent As an example for Theorem \ref{formulaFib}, we calculate $f(4)$ and $f(5)$.
\begin{eqnarray*}
f(4) &=& 1+ \displaystyle \sum_{k=0}^{2} \sum_{i=0}^{2-k} \binom{k}{2-k-i} \\ & =& 1+ \displaystyle \sum_{i=0}^{2} \binom{0}{2-i} + \displaystyle \sum_{i=0}^{1} \binom{1}{1-i} + \displaystyle \sum_{i=0}^{0} \binom{2}{0-i} \\ & =& 1+ \binom{0}{0} + \binom{1}{1} + \binom{1}{0} + \binom{2}{0} \\ & =& 5
\end{eqnarray*}

\vspace*{-0.5cm}
\begin{eqnarray*}
f(5) &=& 1+ \displaystyle \sum_{k=0}^{3} \sum_{i=0}^{3-k} \binom{k}{3-k-i} \\ & =& 1+ \displaystyle \sum_{i=0}^{3} \binom{0}{3-i} + \displaystyle \sum_{i=0}^{2} \binom{1}{2-i} + \displaystyle \sum_{i=0}^{1} \binom{2}{1-i} + \displaystyle \sum_{i=0}^{0} \binom{3}{0-i}\\ & =& 1+ \binom{0}{0} + \binom{1}{1} + \binom{1}{0} + \binom{2}{1} + \binom{2}{0} + \binom{3}{0} \\ & =& 8
\end{eqnarray*}



\section{Conclusion}
In this paper, we give a closed-form representation for Fibonacci numbers using Fibonacci trees. A similar approach can be attempted for finding a closed-form representation for Lucas and Bernoulli numbers.


\section{References}

\begin{thebibliography}{15}


\bibitem{avl} G.M. Adelson-Velskii and E.M. Landis, An algorithm for the organization of information, Soviet Math. Dokl., 3 (1962) 1259-1262.


\bibitem{zeck} R.M. Capocelli, A Note on Fibonacci Trees and the Zeckendorf Representation of Integers, The Fibonacci Quaterly, 26, 4 (1988) 318-324.

\bibitem{Binet} Eric W., ``Binet's Fibonacci Number Formula" from MathWorld.

\bibitem{Optim} Y. Horibe, Notes on Fibonacci Trees and Their Optimality, The Fibonacci Quaterly, 21, 2 (1983) 118-128.

\bibitem{Balaban} N. Jia, K.W. Mclaughlin, Fibonacci trees: A study of the asymptotic behavior of Balaban's index, Communications in Mathematical and in Computer Chemistry (MATCH), 51 (2004) 79-95.

\bibitem{FibFib} S.G. Wagner, The Fibonacci number of Fibonacci trees and a related
family of polynomial recurrence systems, Fibonacci Quart. 45, 3 (2007) 247–253.



%
%
%




\end{thebibliography}
\end{document}